\def\lastrevised{April 2012.}
\documentclass[12pt]{article}
\usepackage{amssymb}

\def\header{Hechler and Laver Trees}
\def\al{\alpha}
\def\anal{{\bf\Sigma}^1_1}
\def\be{\beta}
\def\ff{{\mathcal F}}
\def\hhf{{{\mathbb H}_\ff}}

\def\lam{\lambda}
\def\la{\langle}
\def\llf{{{\mathbb L}_\ff}}
\def\om{\omega}
\def\proof{\par\noindent Proof\par\noindent}
\def\qed{\par\noindent QED\par}
\def\range{{\mbox{ range}}}
\def\rank{{\mbox{rank}}}
\def\ra{\rangle}
\def\res{\upharpoonright}

\def\rmiff{\mbox{ iff }}
\def\si{\sigma}
\def\sm{{\setminus}}
\def\su{\subseteq}

\newtheorem{theorem}{Theorem} 
\newtheorem{lemma}[theorem]{Lemma}
\newtheorem{define}[theorem]{Definition}
\newtheorem{prop}[theorem]{Proposition}
\newtheorem{cor}[theorem]{Corollary}

\begin{document}

\begin{center}
{\large Hechler and Laver Trees}
\end{center}

\begin{flushright}
Arnold W. Miller\footnote{ 
These results  were obtained in December 2002 while
visiting the Fields Institute at the University of Toronto.
Thanks to its director Juris Steprans and his staff
for their hospitality and also for the discussions we had
with Steprans at that time.
It was presented to the SEALS conference in Gainesville
Florida in March 2003, and in a topics course at
the University of Wisconsin, Madison in the fall of 2009.
Last revised \lastrevised.
\par Mathematics Subject Classification 2010: 03E15; 03E05; 03E02
\par Keywords: Analytic sets, trees, constructible sets, Martin's Axiom,
games.}
\end{flushright}

\def\address{\begin{flushleft}
Arnold W. Miller \\
miller@math.wisc.edu \\
http://www.math.wisc.edu/$\sim$miller\\
University of Wisconsin-Madison \\
Department of Mathematics, Van Vleck Hall \\
480 Lincoln Drive \\
Madison, Wisconsin 53706-1388 \\
\end{flushleft}}

\begin{center}
Abstract
\end{center}

\begin{quote}
A Laver tree is a tree in which each node splits infinitely
often.  A Hechler tree is a tree in which each node splits cofinitely
often.   We show that every analytic set is either disjoint from
the branches of a Heckler tree or contains the branches of a 
Laver tree.  As a
corollary
we deduce Silver Theorem that all analytic sets are
Ramsey.  We show that in Godel's constructible universe that
our result is false for co-analytic sets (equivalently it fails for
analytic sets if we switch Hechler and Laver).   We show that
under Martin's axiom that our result holds for $\bf\Sigma_2^1$ sets.
Finally we define two games related to this property.
\end{quote}

\begin{define}
A subtree $H\su \om^{<\om}$ is {\em Hechler} iff
$\forall s\in H\;\forall^\infty n\;\; sn\in H$.
A subtree $L\su \om^{<\om}$ is {\em Laver} iff
$\forall s\in L\;\exists^\infty n\;\;\; sn\in L$.
\end{define}

These definitions are motivated by well-known forcing notions of
Laver \cite{laver} and Hechler \cite{hechler}.  
In the classical Hechler forcing the cofinite sets on the
$n^{th}$ level of the tree would all be the same.

\begin{define}
For any subtree $T\su \om^{<\om}$ define
$$[T]=\{x\in\om^\om\;:\; \forall n\;\; x\res n\in T\}$$
\end{define}

\begin{theorem}\label{eitheror}
For any $\anal$ set $A\su\om^\om$ either there exists a Hechler
tree $H$ with $[H]\cap A=\emptyset$ or there exists a Laver
tree $L$ with $[L]\su A$.
\end{theorem}
\proof
Since analytic sets are projections of closed sets there
exists a tree $T$ on $\om^{<\om}\times\om^{<\om}$ such that
$$A=p[T]=^{def}\{x\in\om^\om\;:
\;\exists y\in\om^\om\; \forall n\; (x\res n,y\res n)\in T\}.$$
Assume that for every Hechler $H$ that $A$ meets $[H]$ and we 
will show there is a Laver $L$ with $[L]\su A$.

For $s,t\in \om^{<\om}$ define
$$A_{s,t}=\{x\in\om^\om\:;\;s\su x \exists y\supseteq t\; (x,y)\in [T]\}.$$

\begin{define}
We say that $H$ is Hechler with root $s$ if
for all $t\in H$ either $s\su t$ or $t\su s$ and beneath
$s$ there is cofinite splitting.
\end{define}

\begin{lemma}\label{lem1}
Suppose for every Hechler $H$ with root $s$ that 
$A_{s,t}\cap [H]\neq\emptyset$.
Then there are infinitely many $n$ such that for every Hechler 
$H$ with root $sn$ that $A_{sn,t}\cap [H]\neq\emptyset$.
\end{lemma}
\proof
Otherwise for all but finitely many $n$ (say $n>N$) there exists
a Hechler $H_n$ with root $sn$ which misses $A_{sn,t}$.  But
then the Hechler tree:
$H=\bigcup_{n>N}H_n$
misses $A_{s,t}$ and has root $s$.
\qed

\begin{lemma}\label{lem2}
Suppose for every Hechler $H$ with root $s$ that $A_{s,t}\cap [H]\neq\emptyset$.
Then there exists an infinite  well-founded tree 
$T\su\{r\;:\; s\su r\}$ with root $s$
and terminal nodes $B\su T$ such that

(1) The nonterminal nodes of $T$ are $\om$-splitting, i.e.,
if $r\in T\sm B$, then there are infinitely many $n$ with $rn\in T$, and

(2) For every $r\in B$ there exists $n$ such that for every Hechler
tree $H$ with root $r$, $A_{r,tn}\cap [H]\neq\emptyset$.
\end{lemma}
\proof

For each ordinal $\al$ define a set $B_\al\su \{r\;:\; s\su r\}$
as follows. 

(a) $r\in B_0$ iff there exists $n$ such that for every Hechler
tree $H$ with root $r$, $A_{r,tn}\cap [H]\neq\emptyset$.

(b) $B_{\al+1}=B_\al\cup\{r\;:\;\exists^\infty n \;\; rn\in B_\al\}$

(c) $B_\lam=\bigcup_{\al<\lam}$ for $\lam$ a limit ordinal.

Define function $\rank(r)$ on $r\supseteq s$ as follows,
$\rank(r)=\al$ if $\al$ is the least ordinal with $r\in B_\al$
and $\rank(r)=\infty$ if there is no such ordinal.

\bigskip\noindent
{\bf Case 1.} $\rank(s)$ is an ordinal.  In this case it is easy
to build $T$ and $B$ as required.

\bigskip\noindent
{\bf Case 2.} $\rank(s)=\infty$.  We show that this is impossible.
Note that if $\rank(r)=\infty$ then for
all but finitely many $n$ we must have that $\rank(rn)=\infty$.
Hence we may construct a Hechler tree $H$ with root $s$ such that
$\rank(r)=\infty$ for all $r\in H$ below the root.  For each $n<\om$ 
for each $r\in \om^{n+|s|}\cap H$ there exists a Hechler $H_r$ with 
root $r$ such that
$$[H_r]\cap A_{s,tn}=\emptyset.$$
Define 
$$K_n=\bigcup\{H_r\;:\; r\in\om^{n+|s|}\cap H\}$$
Note that $K_n$ is a Hechler tree with root $s$ whose
$n+|s|$ level is the same as $H$.  It also true that
$K_n\cap A_{s,tn}=\emptyset$.   Because they are so wide
$K=\bigcap_{n<\om}K_n$ is a Hechler tree with root $s$ such
that $[K]\cap \bigcup_{n}A_{s,tn}=\emptyset$. This contradicts the 
hypothesis of the Lemma since $\bigcup_{n}A_{s,tn}=A_{s,t}$.

\bigskip
Finally, if 
$T$ is trivial, i.e., 
$T=B=\{s\}$
just apply Lemma \ref{lem1} to make $T$ infinite.

\qed

\bigskip\noindent Proof of Theorem \ref{eitheror}:

Suppose for every Hechler tree $H$ with trivial root that
$A\cap [H]\neq\emptyset$.  Apply Lemma \ref{lem2} to obtain
a non-trivial well-found tree $T_0$ with terminal nodes $B_0$
and witnesses of length one.

Suppose we are given a well-founded tree 
$T_n$ with trivial root and terminal nodes $B_n$ such that
for all $s\in T_n\sm B_n$ there are infinitely many immediate
extensions of $s$ in $T_n$ and for each $s\in B_n$ there is
a $t_s$ of length $n+1$ such that for every Hechler
tree $H$ with root $s$, $A_{s,t_s}\cap [H]\neq\emptyset$.
Apply Lemma \ref{lem2} to each node $(s,t_s)$ with $s\in B_n$.
Union all these trees together to get $T_{n+1}$ which end
extends $T_n$.  It follows that $L$ is a Laver where
$$L=\bigcup_{n<\om} T_n$$
Note that although the length of the witnesses grow much
slower than the $s$-part, nevertheless, they union up to show 
that ${L}\su A$.
\qed

\begin{define}
For $\ff$ a filter extending the cofinite filter on $\om$ define
$\hhf$ to be the Hechler trees mod $\ff$, i.e., instead of demanding
that for each $s\in H$ that $sn\in H$ for cofinitely many $n$, we
demand that $$\{n\;:\; sn\in H\}\in\ff.$$  Analogously define
$\llf$ the Laver trees mod $\ff$ by for each $s\in L$ 
$$\{n\;:\; sn\in L\}\in\ff^+$$  where $\ff^+$ are the positive
$\ff$ sets, i.e, 
sets whose complement is not in $\ff$.
\end{define}

\begin{theorem}\label{easy}
For any filter $\ff$ and
any $\anal$ set $A\su\om^\om$ either there exists a Hechler
tree $H\in\hhf$ with $[H]\cap A=\emptyset$ or there exists a Laver
tree $L\in \llf$ with $[L]\su A$.
\end{theorem}
\proof
The proof of this goes over mutatis mutandis, the proof
of Theorem \ref{eitheror}.
\qed

Any Hechler tree $H$ may be pruned so that every
node in it is strictly increasing, i.e. ,
if $\la x_0,x_1,\ldots,x_n\ra\in H$ then
$x_0<x_1<\ldots x_n$.  By the range of $H$ we mean
all infinite subsets of $\om$ which are the image of some branch
$f\in H$, i.e., 
$$\range(H)=\{\{f(n)\;:\;n<\om\}\;:\;f\in [H]\}$$

\begin{prop}\label{proprange}
Suppose $H\in \hhf$.  Then there exists $X\in[\om]^\om$
such that $[X]^\om \su \range(H)$.
\end{prop}
\proof
We may suppose that the nodes of $H$ are strictly increasing.
Construct a strictly sequence $x_0,x_1,\ldots, x_n$
such that for every $k$ and subsequence
$$0\leq i_1<i_2<\cdots<i_k\leq n$$ we have that
$\la x_{i_1},\ldots, x_{i_k}\ra\in H$.  To obtain $x_{n+1}$ we
need only intersect finitely many elements of the filter $\ff$.
\qed

\begin{cor}
(Silver \cite{silver}) Analytic sets have the Ramsey Property.
This means that for any $\anal$ set $A\su[\om]^\om$ there exists
an $X\in [\om]^\om$ with either $[X]^\om\su A$ or 
$[X]^\om\cap A=\emptyset$.
\end{cor}
\proof
Let $\ff$ be a nonprincipal ultrafilter.  Note that
$\hhf=\llf$ for ultrafilters.   Define $B\su\om^\om$ by
$f\in B$ iff $f$ is strictly increasing with range
in $A$.   Then $B$ is $\anal$ and so by Theorem \ref{easy} there is
a Hechler tree $H\in\hhf$ with $[H]\su B$ or $[H]\cap B=\emptyset$.
By Proposition \ref{proprange} there is an infinite $X$ as required.
\qed

This gives a proof of Silver's Theorem which avoids the
accept-reject arguments of Galvin-Prikry \cite{GP} and 
Ellentuck \cite{ellen}.

\begin{theorem}
(V=L) There exists a $\Pi^1_1$ set $A\su\om^\om$ such that
$[H]\cap A\neq\emptyset$ for every Hechler tree $H$ but
$A$ contains no Laver $[L]$.
\end{theorem}
\proof

Using the definable well-ordering of the reals in $L$
construct $B\su\om^\om$ a $\Sigma^1_2$ set with the following
two properties:   

(1) $B$ is an $<^*$ scale, i.e., 
$B=\{g_\al\in\om^\om\;:\; \al<\om_1\}$ where
$\al<\be$ implies $g_\al<^*g_\be$ and for every $f\in\om^\om$ there
exist $\al$ such that $f<^*g_\al$.

(2)  $B$ has the property that for
any $\si:\om^{<\om}\to\om$ there exists $g\in B$ such that
for all $x\in 2^\om$ if $f=2g+x$ then
$\forall n \;\; f(n)>\si(f\res n)$,

Let $C\su \om^\om\times 2^\om$ be $\Pi_1^1$ so that
$g\in B \rmiff \exists x\; (g,x)\in C$.

Given $f\in\om^\om$
define $Q(f)=(g,x)$ where $g\in\om^\om$ and $x\in 2^\om$
are determined by $f=2g+x$.
Define the $\Pi^1_1$ set $A$ by
$$A=\{f\in\om^\om\;:\; Q(f)\in C\}.$$

Note that for any Hechler $H$ we can find $\si:\om^{<\om}\to\om$
such that 
$$H_\si=^{def}\{f\in\om\;:\;\forall n\;\; f(n)>\si(f\res n)\;\}\su H$$
It follows from (2) that $A$ meets
every $[H]$.  On the other hand $A$ cannot contain the branches
$[L]$ of a Laver tree.  This is because of the scale (1).  
Take a 3 splitting subtree of $T\su L$, i.e., for every $s\in T$
there are exactly 3 immediate extensions of $s$ in $T$.
For each $g\in\om^\om$ define
$$C_g=\{f\in\om^\om\;:\;\exists x\in 2^\om\; f=2g+x\}$$
and note that $A\su\bigcup_{g\in B} C_g$.
If $[T]\su A$ then by the scale property of $B$ there would
have to be a countable set $Q\su B$ with such that
$$[T]\su \bigcup_{g\in Q} C_g$$
But the $C_g$ are the branches of a binary splitting tree and
since $T$ is 3-splitting, it is easy to construct $f\in [T]$ such
that $f\notin C_g$ for every $g\in Q$.
\qed

\begin{theorem}
Assume MA +$\neg$CH. If $A\su\om^\om$ is ${\bf \Sigma^1_2}$,
then either there is a Hechler tree $H$ with
$[H]\cap A=\emptyset$ or there is a Laver tree $L$ with
$[L]\su A$.  In fact, this is true for any set $A$ which is 
the union of $\om_1$ many Borel sets.
\end{theorem}
\proof
Suppose $A=\bigcup_{\al<\om_1}B_\al$ where each $B_\al$ is
Borel and the union is increasing.  Since no
$B_\al$ contains the branches of a Laver tree we have
$H_\al$ a Hechler tree with $[H_\al]\cap B_\al=\emptyset$.
Without loss we may assume that 
$$[H_\al]=\{f\in\om^\om\;:\;\forall n \; f(n)>\si_\al(f\res n)\;\}$$
where $\si_\al:\om^{<\om}\to\om$.
By Martin's axiom we may find $\si:\om^{<\om}\to\om$ which
eventually dominates each $\si_\al$.   By a counting argument
we can find a single $\si:\om^{<\om}\to\om$ which everywhere dominates
$\om_1$ of the $\si_\al$.  But this means
that $H_\si$ is a Hechler tree disjoint from $A$ since the
$B_\al$'s are an increasing union.
\qed

\bigskip\bigskip
\noindent Finally we make some remarks about games.

\bigskip\noindent
{\bf Game 1.} Given $A\su\om^\om$.  Player I and II alternatingly
play $$n_0,\; m_0>n_0,\;  \;n_1, \;\; m_1>n_1,\;\; \ldots$$ 
with Player I playing $n_k\in\om$ and
Player II responding with $m_k>n_k$.
The play of the game is won by Player II iff
$(m_i:i<\om)\in A$.

\begin{prop}
(a) Player II has a winning strategy in Game 1 iff there exists
a Laver tree $L$ with $[L]\su A$.
(b) Player I has a winning strategy in Game 1 iff there exists
a Hechler tree $H$ with $[H]\cap A=\emptyset$.
\end{prop}
\proof
Given the trees it easy to get the strategies.
For the other direction:

(a) Use player II's winning strategy to construct
a Laver tree as required.

(b) If $\si:\om^{<\om}\to\om$ is Player I's winning strategy, then
for any sequence $(m_i:i<\om)$ 
such that $m_{i+1}>\si(m_0,\ldots,m_i)$ for every $i$ we have that
$(m_i:i<\om)\notin A$.  But this gives a Hechler tree $H$
with $[H]$ disjoint from $A$.
\qed

\bigskip\bigskip\noindent
{\bf Game 2.}
Given $A\su\om^\om$.  Player I and II alternatingly play
$$X_0,\;\;m_0\in X_1,\;\;X_1,\;\;m_1\in X_1,\ldots$$ 
with Player I playing
$X_k\in[\om]^\om$ and  Player II responding with $m_k\in X_k$.
Player II wins the play of the game iff $(m_i:i<\om)\in A$.

\begin{prop}
(a) Player II has a winning strategy in Game 2 iff there
exists a Hechler tree $H$ with $[H]\su A$.
(b) Player I has a winning strategy in Game 2 iff there exists
a Laver tree $L$ with $[L]\cap A=\emptyset$.
\end{prop}
\proof
From right-to-left in both cases is easy.
For the other direction:

(a) Let $\si$ be a winning strategy for Player II.  
Consider
$$\{m_0\;:\; \exists X_0\;\; \si(X)=m_0\}$$
This set must be cofinite, since otherwise consider
$\si$'s response to its complement.  Similarly given any
sequence $X_0,X_1,\ldots, X_{n-1}$ the set
$$\{m_n\;:\; \exists X_n\;\; \si(X_0,\ldots,X_n)=m_n\}$$
must be cofinite.   Construct $X_s$ for $s\in\om^{<\om}$
and get a Hechler tree $H$ all of whose branches are plays
of the winning strategy and hence are in $A$.

(b) The sequence of $X_s$ played by winning strategy of Player
I determine a Laver tree $L$.

\qed

\bigskip
Some of the results in this note follow from Zapletal \cite{zap}.

\address


\begin{thebibliography}{99}

\bibitem{ellen}
Ellentuck, Erik;
A new proof that analytic sets are Ramsey.
J. Symbolic Logic 39 (1974), 163-165. 

\bibitem{GP} Galvin, Fred; Prikry, Karel;
Borel sets and Ramsey's theorem.
J. Symbolic Logic 38 (1973), 193-198. 

\bibitem{hechler}  Hechler, Stephen H.; On the existence of certain cofinal
subsets of $^{\omega }\omega $. Axiomatic set theory (Proc. Sympos.
Pure Math., Vol. XIII, Part II, Univ. California, Los Angeles, Calif., 1967),
pp. 155-173. Amer. Math. Soc., Providence, R.I., 1974.

\bibitem{laver}
Laver, Richard;
On the consistency of Borel's conjecture.
Acta Math. 137 (1976), no. 3-4, 151-169. 

\bibitem{silver}
Silver, Jack;
Every analytic set is Ramsey.
J. Symbolic Logic 35 (1970), 60-64. 

\bibitem{zap} Zapletal, Jindrich;
Isolating cardinal invariants.
J. Math. Log. 3 (2003), no. 1, 143-162. 

\end{thebibliography}
\end{document}